\documentclass[10pt,onecolumn]{article}

\usepackage{amsmath,amssymb,amsthm}
\newtheorem{theorem}{Theorem}
\newtheorem{lemma}{Lemma}

\newtheorem{cor}{Corollary}

\newcommand{\Nats}{{\mathbb N}}

\newcommand{\E}{{\mathbb E}}
\newcommand{\ER}{\mbox{{Erd\H{o}s-R\'enyi}}}
\newcommand{\prob}{{\mathbb P}}

\newcommand{\wkbar}{{\overline{w^k}}}
\newcommand{\wbar}{{\overline{w}}}
\newcommand{\wtwobar}{{\overline{w^2}}}
\newcommand{\wthreebar}{{\overline{w^3}}}
\newcommand{\wfourbar}{{\overline{w^4}}}
\newcommand{\eps}{\varepsilon}
\newcommand{\Cov}{{\rm Cov}}
\newcommand{\Var}{{\rm Var}}

\newcommand{\degree}{{\rm degree}}

\begin{document}
\title{A model for infection on graphs}

\author{M. Draief\footnote{Department of Electrical and Electronic Engineering
Imperial College London,
Email: M.Draief@imperial.ac.uk}
\and
A. Ganesh\footnote{Department of Mathematics, 
University of Bristol,
Email: A.Ganesh@bristol.ac.uk}
}

\maketitle

{\bf Abstract }
We address the question of understanding the effect of the underlying network topology on the spread of a virus and the dissemination of information when users are mobile performing independent random walks on a graph. To this end we propose a simple model of infection that enables to study the coincidence time of two random walkers on an arbitrary graph. By studying the coincidence time of a susceptible and an infected individual both moving in the graph we obtain estimates of the infection probability.

The main result of this paper is to pinpoint the impact of the network topology on the infection probability. More precisely, we prove that for homogeneous graph including regular graphs and the classical Erd\"os-R\'enyi model, the coincidence time is inversely proportional to the number of nodes in the graph. We then study the model on power-law graphs, that exhibit heterogeneous connectivity patterns, and show the existence of a phase transition for the coincidence time depending on the parameter of the power-law of the degree distribution.

\section{Introduction}

In recent years there have been a surge of hand-held wireless computing devices such as PDAs together with the proliferation of new services. These portable computing devices are equipped with a short-range wireless technology such as WiFi or Bluetooth. Despite providing a great deal of flexibility this ability to wirelessly connect to other devices, and to transfer data on the move, attracted the attention of virus writers who exploit such features for lauching computer-virus outbreaks that take advantage of human mobility \cite{Kleinberg,Leavitt}.

Over the past couple years, there have been indeed reports of malicious code that take advantage of bluetooth vulnerabilities such as the Cabir worm that was detected during the World Athletics Championship \cite{Cabir} and another at a company that has been reported by CommWarrior \cite{CommWarrior}. Despite their small scales, these incidents bode more threats taking advantage of events and locations where individuals gather in close proximity \cite{YanCuellar,SuChanMiklas}.

In much of the literature on mathematical epidemiology, the members of the population
are assumed to occupy fixed locations and the probability of infection passing between a pair of
them in a fixed time interval is taken to be some function of the distance between them. 
Mean-field models are a special case in which this function is a constant  \cite{DG01}. In this work, we
consider a different model in which the agents are mobile and can only infect each other if they
are in sufficiently close proximity. The model is motivated both by certain kinds of biological 
epidemics, whose transmission may be dominated by sites at which individuals gather in close proximity
(e.g. workplaces or public transport for a disease like SARS, cattle markets for foot-and-mouth
disease, etc.) and by malware spreading between wireless devices via Bluetooth connections, for
example.

\noindent 
{\bf Related work.}   Here we briefly describe some of the relevant related
work on modelling epidemic spreading in mobile environments. To our knowledge the first attempts to model virus spreading in mobile networks relies on the use of a non-rigorous mean-field approximations (similar to the classical Kephart-White model \cite{Kephart}) that incorporates the mobility patterns of users. In \cite{Noble}, the authors derive a threshold for the persistence of the epidemic by computing the average number of neighbours of a given  node. Using a similar approach but with different mobility patterns, Nekovee et al. \cite{Nekovee,NekoveeRhodes}  explore the evolution of the number of devices that are infected in terms of the contact rate betwen users. 

A related line of work studying the dissemination of information in opportunistic networks \cite{Chaintreau} focuses on the following analogous problem:  Suppose that all individuals are interested in a piece of information that is initially held by one user. The information is transmitted between users who happen to be close to each other. As in the case of static networks \cite{Pittel}, one may be interested in the time it takes for the rumour to be known to all users. To this end we need to understand how information is transmitted between an informed and an ignorant user. Our work gives some insight on the impact of the network structure on the likelihood of successfully transmitting the rumour.

\noindent 
{\bf Our contribution.} In contrast to the previous work which has focused on Euclidean models and homogeneous mobility patterns, in this work we consider a model wherein the different locations that a user can reach have varying popularity.

More precisely, we consider a simple and stylised mathematical model of the spread of infection as follows. There
is a finite, connected, undirected graph $G=(V,E)$ on which the individuals perform independent 
random walks: they stay at which vertex for an exponentially distributed time with unit mean, 
and then move to a neighbour of that vertex chosen uniformly at random. The infection can pass 
from an infected to a susceptible individual only if they are both at the same vertex, and the 
probability of its being passed over a time interval of length $\tau$ is $1-\exp(-\beta \tau)$, 
where $\beta>0$ is a parameter called the infection rate. We shall consider a single infected 
and a single susceptible individual and ask what the probability is that the susceptible 
individual becomes infected by time $t$. This probability has been studied in the case of 
a complete graph in \cite{DD04}. Here, we extend their results to a much wider class of graphs.

It is simplistic to consider just a single infective and a single susceptible individual.
Nevertheless, insights gained from this setting are relevant in the ``sparse" case, where 
the number of both infected  and susceptible individuals is small and inter-contact times 
are fairly large. In that case, it is not a bad approximation to consider each pair of
individuals in isolation. The ``dense" setting will require quite different techniques 
and is not treated here.

The rest of the paper is organised as follows. In Section \ref{se-model}, we present our model and the family of networks we will consider. Besides we state our main results that relate the coincidence time of the two walkers to the stationary distribution of a random walk on a graph. In Section \ref{se-results} we give a detailed proof of our main result on the probability of infection for regular graphs, the Erd\"os-R\'enyi graph and power-law networks. Section \ref{se-conclusion}
summarises our contribution and suggests further extensions of our work.

\section{Models and results}\label{se-model}

We now describe the model precisely.
Let $X_t, Y_t \in V$ denote the positions of the susceptible and infected individuals
respectively at time $t$. We model $(X_t, t\ge 0)$ and $(Y_t, t\ge 0)$ as independent 
continuous-time Markov chains (CTMCs) on the finite state space $V$, with the same 
transition rate matrix given by
$$
q_{xy} = \left\{
\begin{array}{ll}
\frac{1}{\degree(x)} & \mbox{if $(x,y) \in E$}, \\
0 & \mbox{if $y\neq x$ and $(x,y) \notin E$}, \\
-1 & \mbox{if $y=x$}
\end{array}
\right.
$$

where ${\rm degree}(x)$ is the number of neighbours of $x$ (nodes $y$ such that $(x,y)\in E$) in the graph $G$.

We define the coincidence time up to time $t$, denoted $\tau(t)$, as the total
time up to $t$ during which both walkers are at the same vertex, i.e.,
\begin{equation} \label{def_coincide}
\tau(t) = \int_0^t 1(X_s=Y_s)ds.
\end{equation}
Let $\gamma(t)$ denote the probability that the initial susceptible becomes infected 
by time $t$. Then, conditional on $\tau(t)$, we have
\begin{equation} \label{infect_prob}
\gamma(t) = 1-\exp(-\beta \tau(t)),
\end{equation}
where $\beta>0$ is the infection rate.

We are interested in estimating the coincidence time $\tau(t)$ and the
infection probability $\gamma(t)$ for different families of graphs.

Observe that the Markov chains $X_t$, $Y_t$ have invariant distribution 
$\pi$ given by 
\begin{equation} \label{invar_dist}
\pi_x= \frac{\degree(x)}{\sum_{v\in V} \degree(v)}
\end{equation} 
and that they are reversible, i.e., $\pi_x q_{xy} = \pi_y q_{yx}$ 
for all $x,y \in V$.

We consider the case when these chains are started independently
in the stationary distribution and provide estimates on the coincidence
time and the infection probability, for arbitrary graphs. 

\begin{theorem} \label{thm:basic}
Suppose $X_0$ and $Y_0$ are chosen independently according to the
invariant distribution $\pi$. Then, we have
$$
\E [\tau(t)] = \sum_{v\in V} \pi_v^2 t, \mbox{ and }
\E [\gamma(t)] \leq 1-\exp \Bigl( -\beta t\sum_{v\in V} \pi_v^2 \Bigr).
$$
\end{theorem}

\noindent \emph{Proof.}
Observe that, for all $s\ge 0$,
$$
\prob(X_s=Y_s) = \sum_{v\in V} \prob(X_s=Y_s=v) = \sum_{v\in V} \pi_v^2,
$$
because $X_s$ and $Y_s$ are independent, and are in stationarity.
Hence, it is immediate from (\ref{def_coincide}) that
$$
\E[\tau(t)] = \int_0^t \sum_{v\in V} \pi_v^2 ds = \sum_{v\in V} \pi_v^2 t.
$$
Next, taking expectations in (\ref{infect_prob}) with respect to the
conditioning random variable $\tau(t)$, we have 
\begin{eqnarray*}
\E [\gamma(t)] &=& 1 - \E [\exp(-\beta \tau(t))]\\
&\leq & 1 - \exp(-\beta \E [\tau(t)]) \\
&=& 1-\exp \Bigl( -\beta t\sum_{v\in V} \pi_v^2 \Bigr),
\end{eqnarray*}
where the inequality follows from Jensen's inequality.
\hfill $\Box$

We now introduce some terminology
and define some examples of graph models that we shall consider.

For two functions $f(\cdot)$ and $g(\cdot)$ on the natural numbers, we write 
$f(n) \sim g(n)$ to mean that their ratio tends to 1 as $n$ tends to infinity. 
We write $f(n)=O(g(n))$ if $f(n)/g(n)$ remains bounded by a finite constant,
$f(n)=o(g(n))$ if $f(n)/g(n)$ tends to zero, and $f(n)=\Omega(g(n))$ if
$g(n)=O(f(n))$. 
For a sequence of events $A_n$ indexed by $n\in \Nats$, we  say that they occur 
with high probability (whp) if $\prob(A_n)$ tends to 1 as $n$ tends to infinity. 

\medskip

\noindent{\bf Examples}

\noindent
\emph{Complete graphs. }  
Consider the complete graph on $n$ nodes, namely the graph in
which there is an edge between every pair of nodes. Thus, $\degree(v)=n-1$ and
$\pi_v=1/n$ for all $v\in V$, so we have by Theorem \ref{thm:basic} that
$\E[\tau(t)] = t/n$. This result should be intuitive by symmetry. Lemma \ref{thm:basic}
also gives us an upper bound on the infection probability, $\E[\gamma(t)] \leq
1-\exp(-\beta t/n)$. Roughly speaking, this says that it takes time of order $n/\beta$
for the susceptible individual to become infected; for $t\ll n/\beta$, the probability
of being infected is vanishingly small. Again, this is consistent with intuition.

\medskip

\noindent
\emph{Regular graphs} A graph $G=(V,E)$ is said to be $r$-regular if $\degree(v)=r$
for all $v\in V$. Thus, a complete graph is regular with $r=n-1$. 
It is readily verified that $\pi_v=1/n$ for all $v\in V$ if $G$ is
$r$-regular for any $r\ge 2$. (If $r=1$, then $G$ is a matching and is not connected.)
Hence, if $G$ is connected, we have the same estimates for $\tau(t)$ and $\gamma(t)$
as for the complete graph, which is a special case corresponding to $r=n-1$.

The next examples we consider will be families of random graphs widely used in practice
to model networks.
\medskip

\noindent
\emph{$\ER$ random graphs} The $\ER$ graph $G(n,p)$ is defined as a random graph on 
$n$ nodes, wherein each edge is present with probability $p$, independent of all other 
edges. We consider a family of such random graphs indexed by $n$, and take $p$ to be a 
function of $n$ chosen so that $np > c\log n$ for some constant $c>1$. We also condition 
on the graph being connected. For $p$ as above, the probability of connectivity tends to 
1 as $n$ tends to infinity, so conditioning on connectivity does not alter any of the 
estimates we shall derive later for the coincidence time on such graphs. In this model, 
the node degrees are identically distributed Binomial random variables with parameters 
$(n-1,p)$. In particular, they concentrate around the mean value of $np$, and have
exponentially decaying tails away from this value. Thus, while $\ER$ graphs are not 
exactly regular, they exhibit considerable homogeneity in node degrees. 

\medskip

\noindent
\emph{Power law random graphs} In contrast to the above graph models, many real-world 
networks exhibit considerable heterogeneity in node degrees, and have empirical degree 
distributions whose tails decay polynomially; see, e.g., \cite{BA,FFF}. This observation 
has led to the development of generative models for graphs with power-law tails \cite{BA,BR} 
as well as random-graph models possessing this property \cite{ChLu03}. For definiteness, 
we work with the model proposed in \cite{ChLu03}, but we believe that similar results will
hold for the other models as well.

In the model of \cite{ChLu03}, each node $v$ is associated with a positive 
weight $w_v$, and edges are present independently with probabilities related 
to the weights by
\begin{equation} \label{edgeprob}
\prob((u,v)\in E) = \frac{w_u w_v}{W} \mbox{ where } W=\sum_{x\in V} w_x.
\end{equation}
We assume that $W\geq w_{\max}^2$, so that the above defines a probability.
It can be verified that $\E [\degree(v)]=w_v$ and so
this model is also referred to as the expected degree model. The model allows 
self-loops. The $\ER$ graph $G(n,p)$ is a special case corresponding to the
choice $w_v=np$ for all $v\in V$.
If the weights are chosen to have a power-law distribution, then 
so will the node degrees. The following 3-parameter model for the ordered weight
sequence is proposed in \cite{ChLu03}, parametrised by the mean degree $d$, the
maximum degree $m$, and the exponent $\gamma>2$ of the weight distribution:
\begin{equation} \label{weights}
w_i = m \Bigl( 1+\frac{i}{i_0} \Bigr)^{-\frac{1}{\gamma-1}}, \quad i=0,1,\ldots,n-1,
\end{equation}
where
\begin{equation} \label{i0}
i_0 = n \Bigl( \frac{d(\gamma-2)}{m(\gamma-1)} \Bigr)^{\gamma-1}.
\end{equation}
Note that $W=\sum_{i=0}^{n-1} w_i \sim nd$.

We consider a sequence of such graphs indexed by $n$. The maximum expected degree $m$ 
and the average expected degree $d$ may, and indeed typically will, depend on $n$. In
models of real networks, we can typically expect $d$ to remain bounded or to grow slowly
with $n$, say logarithmically, while $m$ grows more quickly, say as some fractional power
of $n$. In this paper, we only assume the following:
\begin{equation} \label{assumptions}
d \ge \delta > 0, \; d = o(m), \; m \leq \sqrt{nd}, 
\; \frac{m}{d} = o\Bigl( n^{\frac{1}{\gamma-1}} \Bigr).
\end{equation}
Here, $\delta$ is a constant that does not depend on $n$. In other words, the average 
expected degree is uniformly bounded away from zero. The third assumption simply
restates the requirement that $w_0^2 \leq W$, so that (\ref{edgeprob}) defines valid
probabilities. The last assumption ensures that $i_0$, defined in (\ref{i0}), tends
to infinity.

We now describe our results about these models.

\begin{theorem} \label{thm:graphs}
Consider a sequence of graphs $G=(V,E)$ indexed by $n=|V|$. On each graph, consider
two independent random walks with initial condition $X_0$, $Y_0$ chosen independently
from the invariant distribution $\pi$ for the random walk on that graph.

We have $\E [\tau(t)] = t/n$ for regular graphs, including the complete graph,
on $n$ nodes. 

For $\ER$ random graphs $G(n,p)$ conditioned to be connected,
and having $np\geq c\log n$ for some $c>1$, we have $\E [\tau(t)] \sim t/n$,
as $n$ tends to infinity. 

Finally, consider a sequence of power law random graphs 
defined via (\ref{edgeprob}) and (\ref{weights}), and satisfying the assumptions
in (\ref{assumptions}). Then, we have the following:
$$
\frac{n \E [\tau(t)]}{t} \sim \begin{cases}
c, & \mbox{if } \gamma>3, \cr
c(\log m), & \mbox{if } \gamma=3, \cr
c(m.d)^{3-\gamma}, & \mbox{if } 2<\gamma<3,
\end{cases}
$$
where $c>0$ is a constant that may depend on $\gamma$, but not on $n$, $m$ or $d$.
\end{theorem}

\section{Proof of Theorem \ref{thm:graphs}}\label{se-results}

If the graph $G$ is regular, then, by (\ref{invar_dist}), $\pi_v=1/n$ for all
$v\in V$. Hence, the claim of the theorem follows from Theorem \ref{thm:basic}.

In order to estimate $\E[\tau(t)]$, we need to compute 
\begin{equation} \label{probssq}
\sum_{v\in V} \pi_v^2 = 
\frac{\sum_{v\in V} \degree(v)^2}{\Bigl( \sum_{v\in V}\degree(v) \Bigr)^2}.
\end{equation}
Define
\begin{equation} \label{deg_sum_def}
D = \sum_{v \in V} \degree(v) = \sum_{u,v\in V} A_{uv}, 
\end{equation}
where  $A_{uv}=1((u,v)\in E)$, and 

\begin{eqnarray} \nonumber
X_v &=& \degree(v)(\degree(v)-1) = \sum_{i\neq j} A_{vi} A_{vj}\\ \label{d2def}
 D_2 &=& \sum_{v\in V} X_v.
\end{eqnarray}

We will derive the first and second moments of the variables $D$ and $D_2$. It then suffices to use Chebyshev's inequality to establish concentration results for both variables $D$ and $D_2$. By
(\ref{probssq}) and Theorem \ref{thm:basic}, and the fact that $\sum_{v\in V} \degree(v)^2 = D_2+D$ and $D = \sum_{v \in V} \degree(v)$, we will have an estimate of the coincidence time
that holds whp.

We begin by computing the mean and variance of $D$ in the expected degree
model with arbitrary weight sequence $\{ w_i, i=0,\ldots,n-1 \}$.

For notational convenience, we define
$$
\wbar=\frac{1}{n}\sum_{i=0}^{n-1} w_i, \quad
\wkbar = \frac{1}{n} \sum_{i=0}^{n-1} w_i^k, \; k=2,3,\ldots
$$
We obtain $\ER$ graphs $G(n,p)$ by setting $w_i=np$ for all $i$, and so,
$\wkbar = (np)^k$
for such graphs. 

Next, consider power-law graphs with weight sequence specified by (\ref{weights}) and 
(\ref{i0}). Since $i_0$ tends to infinity by assumption, we have for such graphs that
\begin{eqnarray}
\wkbar &=& 
\frac{m^k}{n} \sum_{i=0}^{n-1} \Bigl( 1+\frac{i}{i_0} \Bigr)^{-\frac{k}{\gamma-1}} \nonumber \\
&\sim& \frac{m^k}{n} \int_{0}^{n} \Bigl( 1+\frac{x}{i_0} \Bigr)^{-\frac{k}{\gamma-1}} dx \nonumber \\
&=& m^k \frac{i_0}{n} \int_0^{n/i_0} (1+x)^{-\frac{k}{\gamma-1}} dx. \label{wkbar}
\end{eqnarray}
Now, straightforward calculations yield that $\wbar \sim d$ for all $\gamma>2$, whereas,
for $k\ge 2$, we have
\begin{equation} \label{wkbar_asymp}
\wkbar \sim \begin{cases}
\frac{(\gamma-2)^{k}}{(\gamma-1)^{k-1}(\gamma-1-k)} d^k, & \mbox{if } \gamma>k+1, \cr
\frac{(k-1)^{k}}{k^{k-1}} d^k \log \frac{m}{d}, & \mbox{if } \gamma=k+1, \cr
\frac{(\gamma-2)^{\gamma-1}}{(\gamma-1)^{\gamma-2}(k+1-\gamma)} 
d^{\gamma-1} m^{k+1-\gamma}, & \mbox{if } 2<\gamma<k+1.
\end{cases}
\end{equation}
We can now compute the mean and variance of $D$, the sum of node degrees.

\begin{lemma} \label{lem:total_deg_stats}
Consider a random graph $G=(V,E)$ specified by the expected degree model with
an arbitrary weight sequence $\{ w_v, v\in V \}$ satisfying $W\geq w_{\max}^2$,
where $W=\sum_{v\in V} w_v$. Let the sum of node degrees, $D$, be defined as in
(\ref{deg_sum_def}). Then, we have
\begin{eqnarray} \nonumber
\E[D] &=& n\wbar, \\ \label{eq:total_deg_stats} 
\Var(D) &=& 2\Bigl( n\wbar - \Bigl( \frac{\wtwobar}{\wbar} \Bigr)^2 \Bigr)
- \Bigl( \frac{\wtwobar}{\wbar} - \frac{\wfourbar}{n\wbar^2} \Bigr),
\end{eqnarray}
where $n=|V|$ is the total number of nodes. 

In particular, if $G$ is the
$\ER$ random graph $G(n,p)$, then
\begin{eqnarray} \nonumber 
\E[D] &=& n^2 p\\ \label{eq:total_degstats_er}
 \Var(D) &=& (2n-1)np(1-p) \sim 2n^2 p(1-p),
\end{eqnarray}

whereas, if $G$ is a power law random graph satisfying the assumptions
of Theorem \ref{thm:graphs}, then 

$$\E[D] = nd \quad \text{while} \quad \Var(D) \sim 2nd, \: \text{whp}.$$
.
\end{lemma}

\noindent{\em Proof.}
It is immediate from (\ref{deg_sum_def}) that
\begin{eqnarray*}
\E[D]& = &\sum_{u,v\in V} \prob((u,v)\in E)\\& =& \sum_{u,v\in V} \frac{w_u w_v}{W} = W,
\end{eqnarray*}
which establishes the first equality in (\ref{eq:total_deg_stats}). 
Next, rewrite (\ref{deg_sum_def}) as
$$
D= 2 \sum_{i=1}^n \sum_{j=i+1}^n A_{ij} + \sum_{i=1}^n A_{ii},
$$
and observe from the independence of the edges that
\begin{eqnarray*}
\Var(D) &=& 4\sum_{i=1}^n \sum_{j=i+1}^n \Var(A_{ij}) + \sum_{i=1}^n \Var(A_{ii}) \\
&=& 2\sum_{i=1}^n \sum_{j=1}^n \Var(A_{ij}) - \sum_{i=1}^n \Var(A_{ii}). 
\end{eqnarray*}
Now, $\Var(A_{uv}) = \prob((u,v)\in E)(1-\prob((u,v)\in E))$, and so,
$$
\Var(D) = 2\sum_{i=1}^n \sum_{j=1}^n \Bigl( \frac{w_i w_j}{W} - 
\frac{w_i^2 w_j^2}{W^2} \Bigr) - \sum_{i=1}^n \Bigl(
\frac{w_i^2}{W} - \frac{w_i^4}{W^2} \Bigr).
$$
Upon simplifying, this yields the second equality in (\ref{eq:total_deg_stats}).
Now, using the fact that $\wkbar = (np)^k$ for $\ER$ graphs $G(n,p)$, we readily
obtain (\ref{eq:total_degstats_er}).

Next, suppose $G$ is a power-law graph (more precisely, $G_n$ is a sequence of
power law graphs) satisfying the assumptions of Theorem \ref{thm:graphs}. It follows
from (\ref{wkbar}) that
\begin{equation} \label{w4bar1}
\wfourbar \sim \begin{cases}
\frac{(\gamma-2)^4}{(\gamma-1)^3(\gamma-5)} d^4, & \mbox{if } \gamma>5, \cr
\frac{81}{64} d^4 \log \frac{m}{d}, & \mbox{if } \gamma=5, \cr
\frac{(\gamma-2)^{\gamma-1}}{(\gamma-1)^{\gamma-2}(5-\gamma)} d^{\gamma-1} m^{5-\gamma},
& \mbox{if } 2<\gamma<5,
\end{cases}
\end{equation}
while
\begin{equation} \label{w2bar1}
\wtwobar \sim \begin{cases}
\frac{(\gamma-2)^2}{(\gamma-1)(\gamma-3)} d^2, & \mbox{if } \gamma>3, \cr
\frac{1}{2} d^2 \log \frac{m}{d}, & \mbox{if } \gamma=3, \cr
\frac{(\gamma-2)^{\gamma-1}}{(\gamma-1)^{\gamma-2}(3-\gamma)} d^{\gamma-1} m^{3-\gamma},
& \mbox{if } 2<\gamma<3,
\end{cases}
\end{equation}
and $\wbar \sim d$ for all $\gamma>2$. 

By (\ref{eq:total_deg_stats}), it suffices to show that 
\begin{equation} \label{eq:moment_bounds1}
\Bigl( \frac{\wtwobar}{\wbar} \Bigr)^2 = o(nd) \mbox{ and }
\frac{\wfourbar}{\wbar^2} =o(n^2 d) 
\end{equation}
in order to show that $\Var(D) \sim 2n\wbar \sim 2nd$. 

Suppose first that $\gamma > 3$. Then, by (\ref{w2bar1}) and the fact that $\wbar=nd$, 
$$
\frac{1}{nd} \Bigl( \frac{\wtwobar}{\wbar} \Bigr)^2 = O \Bigl( \frac{d}{n} \Bigr)= o(1),
$$
where the last equality follows by (\ref{assumptions}), and the fac that $d\leq n$.

Now let $\gamma = 3$. Then, by (\ref{w2bar1}) and the fact that $\wbar=nd$, 
$$
\frac{1}{nd} \Bigl( \frac{\wtwobar}{\wbar} \Bigr)^2 = O \Bigl( \frac{d}{n}\log \frac{m}{d} \Bigr)
= O \Bigl( \frac{m}{n} \frac{d}{m}\log \frac{m}{d} \Bigr) = o(1),
$$
where the last equality follows by (\ref{assumptions}). On the other hand, if $2<\gamma<3$, then,
by (\ref{w2bar1}),
\begin{eqnarray*}
\frac{1}{nd} \Bigl( \frac{\wtwobar}{\wbar} \Bigr)^2 &=&
O \Bigl( \frac{d^{2\gamma-5}m^{6-2\gamma}}{n} \Bigr) \\ &=&
O \Bigl( \Bigl( \frac{d}{n} \Bigr)^{\gamma-2} \Bigr) = o(1),
\end{eqnarray*}
where we have used the inequality $m\leq \sqrt{nd}$ from (\ref{assumptions}) to obtain
the second equality. To obtain the last equality, note that it follows from (\ref{assumptions}) 
that $m=o(n)$ and hence that $d=o(n)$ as well. We have thus established the first equality
in (\ref{eq:moment_bounds1}) for all $\gamma>2$. The proof of the second equality is similar
and is omitted. This completes the proof of the lemma.
\hfill $\Box$

The following corollary is now an easy consequence of Chebyshev's inequality.

\begin{cor}
If $G_n$, $n\in \Nats$ is a sequence either of $\ER$ random graphs or of power-laws random
graphs satisfying the assumptions of Theorem \ref{thm:graphs}, then the sum of node degrees
$D$ concentrates at its expected value in the sense that $D\sim \E[D]$ with high probability (whp).
\end{cor}

We now establish a similar concentration result for the sum of squared degrees. To this end, recall that
\begin{eqnarray*}
X_v &=& \degree(v)(\degree(v)-1) = \sum_{i\neq j} A_{vi} A_{vj}\\ 
 D_2 &=& \sum_{v\in V} X_v.
\end{eqnarray*}
 We have the following:

\begin{lemma} \label{lem:d2stats}
Let $D_2$ be defined as in (\ref{d2def}). We then have
\begin{eqnarray} \nonumber
\E[D_2] &=& n\wtwobar - \Bigl( \frac{\wtwobar}{\wbar} \Bigr)^2\\ \label{eq:d2stats}
\Var(D_2) &\le& 4n\wthreebar + 2n\wtwobar + 4n \frac{ (\wtwobar)^2 }{\wbar}.
\end{eqnarray}
\end{lemma}

\noindent \emph{Proof.}
We first note that
\begin{eqnarray*}
\E[X_v] &=& \sum_{i\neq j} \frac{w_i w_j w_v^2}{W^2}\\ & = &
w_v^2 \Bigl( 1-\frac{1}{W^2} \sum_{i\in V} w_i^2 \Bigr) \\ &= &
w_v^2 \Bigl( 1-\frac{\wtwobar}{n\wbar^2} \Bigr).
\end{eqnarray*}
Therefore,
$$
\E[D_2] = \sum_{v\in V} \E[X_v] = n\wtwobar - \Bigl( \frac{\wtwobar}{\wbar} \Bigr)^2,
$$
which is the first part of (\ref{eq:d2stats}).
Next, for distinct nodes $u,v \in V$, we have
\begin{eqnarray*}
\Cov(X_u,X_v) &=& \sum_{i\neq j} \sum_{k\neq l} 
\Cov(A_{iu}A_{ju},A_{kv}A_{lv}) \\
&=& 4 \sum_{i\neq v, l\neq u} 
\Cov(A_{iu}A_{uv},A_{uv} A_{lv}) \\
&=& 4 \E[A_{u,v}] (1-\E[A_{u,v}])  \sum_{i\neq v} \E[A_{iu}] \sum_{l\neq u} \E[A_{lv}].
\end{eqnarray*}
The second equality above holds because, by the independence of edges, the indicator
random variables $A_{iu} A_{ju}$ and $A_{kv} A_{lv}$ corresponding to the open triangles
(or 2-stars) $iuj$ and $kvl$ are independent unless two of the edges are the same; the
only way this can happen is if $(u,v)$ is a common edge and there are 4 possible node
labellings corresponding to each such edge set. Now, recall that $\E[A_{u,v}]=w_u w_v/W$
and $\sum_i \E[A_{iu}] = \E[\degree(u)] = w_u$. Hence, we see from the above that
\begin{equation} \label{covx1}
0 \le \Cov(X_u,X_v) \le 4 \frac{w_u^2 w_v^2}{W}.
\end{equation}
Similarly, we obtain
\begin{eqnarray*}
\Var(X_u) &=& \sum_{i\neq j} \sum_{k\neq l} \Cov(A_{iu}A_{ju},A_{ku}A_{lu}) \\
&=& 4 \sum_j \sum_{i\neq j} \sum_{l\neq i,j} \Cov(A_{iu}A_{ju},A_{ju} A_{lu}) \\&&
+ 2 \sum_{i\neq j} \Var(A_{iu} A_{ju}) \\
&\le& 4 \sum_j \sum_{i\neq j} \sum_{l\neq i,j} \E[A_{iu} A_{ju} A_{lu}]\\&&
+ 2 \sum_{i\neq j} \E[A_{iu} A_{ju}].
\end{eqnarray*}
Using the fact that distinct edges are independent, we get
\begin{equation} \label{varx1}
\Var(X_u) \le 4w_u^3 + 2w_u^2.
\end{equation}
Now, by (\ref{d2def}), (\ref{covx1}) and (\ref{varx1})
\begin{eqnarray*}
\Var(D_2) &=& \sum_{u\in V} \Var(X_u) + \sum_{u\neq v} \Cov(X_u,X_v) \\
&\le& \sum_{u\in V} (4w_u^3+2w_u^2) + \sum_{u,v \in V} 4\frac{w_u^2 w_v^2}{W}.
\end{eqnarray*}
Computing the above sums yields the second part of (\ref{eq:d2stats}).
\hfill $\Box$

We now specialise the results to $\ER$ and power law random graphs, showing that
$D_2$ concentrates near its expected value with high probability.

\begin{lemma} \label{lem:twostar_er}
Suppose $G(n,p)$ is a sequence of $\ER$ random graphs indexed by $n$ (where $p$ depends on 
$n$ but this is not made explicit in the notation), and that $np$ is uniformly bounded away 
from zero. Then $D_2 \sim \E[D_2] \sim n^3 p^2$ whp.
\end{lemma}

\noindent \emph{Proof.}
We have, by Lemma \ref{lem:d2stats} and the fact that $\wkbar = (np)^k$ for the
$\ER$ random graph $G(n,p)$, that
$$
\E[D_2] = n^2(n-1)p^2 \sim n^3 p^2, \quad \Var(D_2) \le 8n^4 p^3 + 2n^3 p^2.
$$
Hence, by Chebyshev's bound, we obtain for all $\eps>0$ that,
\begin{eqnarray*}
\prob ( |D_2-\E[D_2]|>\eps \E[D_2] ) &\le& \frac{\Var(D_2)}{\eps^2 \E[D_2]^2}\\ &\le&
\frac{1}{\eps^2 (n-1)^2 p} \\&&+ \frac{1}{\eps^2 n(n-1)^2 p^2}.
\end{eqnarray*}
Now, by the assumption that $np$ is bounded away from zero, $(n-1)^2 p$ and
$n(n-1)^2 p^2$ tend to infinity as $n$ tends to infinity. Thus, 
$\prob ( |D_2-\E[D_2]|>\eps \E[D_2] )$ tends to zero for all $\eps>0$.
This establishes the claim of the lemma.
\hfill $\Box$

\begin{lemma} \label{lem:twostar_plg}
Suppose $G_n, n\in \Nats$ is a sequence of random graphs satisfying the assumptions
in Theorem \ref{thm:graphs}, with $\gamma >2 $. Then, $D_2 \sim \E[D_2]$ whp, and
$$
\E[D_2] \sim \begin{cases}
cnd^2, & \mbox{if } \gamma>3, \cr
cnd^2(\log m), & \mbox{if } \gamma=3, \cr
cn d^{\gamma-1} m^{3-\gamma}, & \mbox{if } 2<\gamma<3,
\end{cases}
$$
\end{lemma}

\noindent \emph{Proof.}
We will show that $\Var(D_2)=o(\E[D_2]^2)$, so that the claim follows by Chebyshev's
bound, as in the proof of the previous lemma. We will consider separately the
parameter ranges $\gamma \ge 4$, $3\le \gamma < 4$ and $2<\gamma<3$, where $\gamma$
is the exponent in the power law describing the degree distribution.

In the following, $c_1, c_2, \ldots$ will denote generic positive constants, not
necessarily the same from line to line. Recall that $\wbar \sim d$.

Suppose first that $\gamma \ge 4$. Then, by (\ref{wkbar_asymp}), 
$\wthreebar = O(d^3 \log \frac{m}{d})$ and $\wtwobar \sim c_1 d^2$. 
Therefore, by Lemma \ref{lem:d2stats},

$$
\E[D_2] \sim c_1 nd^2$$
and $$ \Var(D_2)=O \Bigl( nd^3\log \frac{m}{d}+nd^2 \Bigr)
=O \Bigl( nd^3 \log \frac{m}{d} \Bigr),
$$
where the last equality holds because of the assumption in (\ref{assumptions})
that $d\ge \delta$ for some constant $\delta>0$. Thus, we see that
$$
\frac{\Var(D_2)}{\E[D_2]^2} = O \Bigl( \frac{1}{nd}\log \frac{m}{d} \Bigr) = o(1),
$$
since $m\le n$.

Suppose next that $3\le \gamma <4$. Then, by (\ref{wkbar_asymp}), $\wthreebar =
O(d^{\gamma-1} m^{4-\gamma})$, while $\wtwobar \sim c_1 d^2$ if $3<\gamma<4$
and $\wtwobar \sim c_2 d^2 \log \frac{m}{d}$ if $\gamma=3$. Therefore, by
Lemma \ref{lem:d2stats},
\begin{eqnarray} \nonumber 
\E[D_2] & \ge & c_1 nd^2 - c_2 d^2 \log^2 \frac{m}{d} \\ \label{ed2_3gamma4} &\ge & c_1 nd^2 - c_2 d^2 \log^2 n
\; = \; \Omega(nd^2),
\end{eqnarray}
whereas,
\begin{eqnarray*}
\Var(D_2) &\le& c_1 nd^{\gamma-1} m^{4-\gamma} + c_2 nd^2 \log \frac{m}{d} 
+ c_3 nd^3 \log^2 \frac{m}{d} \\
&\le& c_1 nd^{\gamma-1} m^{4-\gamma} + c_2 nd^3 \log^2 \frac{m}{d}, \\
&=& c_1 nd^{\gamma-1} m^{4-\gamma} \Bigl( 1+ \Bigl( \frac{d}{m} \Bigr)^{4-\gamma}
\log^2 \frac{m}{d} \Bigr).
\end{eqnarray*}
We have used the assumption that $d$ is uniformly bounded away from zero to obtain 
the second inequality above. Since we also assumed in (\ref{assumptions}) that
$d=o(m)$, we have $(d/m)^{4-\gamma} \log^2(m/d) = o(1)$ for all $\gamma<4$. Hence,
$\Var(D_2)=O(nd^{\gamma-1}m^{4-\gamma})$. Combining this with (\ref{ed2_3gamma4}), we get
$$
\frac{\Var(D_2)}{\E[D_2]^2} =  O \Bigl( \frac{1}{nd} \Bigl( \frac{m}{d} \Bigr)^{4-\gamma} \Bigr)
= O \bigl(\frac{1}{nd} n^{(4-\gamma)/(\gamma-1)} \bigr) = o(1).
$$
We have used (\ref{assumptions}) to obtain the second equality above and the fact that
$\gamma \ge 3$ to obtain the last equality. Moreover, $\E [D_2]\sim c nd^2$ for $3<\gamma<4$, whereas $\E [D_2]\sim c nd^2\log(m)$ for $\gamma=3$.

Finally, suppose that $2<\gamma<3$. Then, by (\ref{wkbar_asymp}), $\wthreebar =
O(d^{\gamma-1} m^{4-\gamma})$ and $\wtwobar \sim c_1 d^{\gamma-1} m^{3-\gamma}$,
so that, by Lemma \ref{lem:d2stats},
\begin{eqnarray*}
\E[D_2] &\ge & c_1 nd^{\gamma-1}m^{3-\gamma} - c_2 (d^{\gamma-2}m^{3-\gamma})^2\\
& \ge & c_1 nd^{\gamma-1}m^{3-\gamma} \Bigl( 1-\frac{c_2}{n} \Bigl( \frac{m}{d} \Bigr)^{3-\gamma} \Bigr).
\end{eqnarray*}
Now, by (\ref{assumptions}), $(m/d)^{3-\gamma} = o(n^{(3-\gamma)/(\gamma-1)}) = o(n)$ since
$\gamma>2$. Consequently,
$$
\E[D_2] = \Omega( nd^{\gamma-1}m^{3-\gamma} ).
$$
On the other hand,
\begin{eqnarray*}
\Var(D_2) &\le& c_1 nd^{\gamma-1} m^{4-\gamma} + c_2 nd^{\gamma-1}m^{3-\gamma}\\
&& + c_3 nd^{2\gamma-3}m^{6-2\gamma} \\
&\le& c_1 nd^{\gamma-1} m^{4-\gamma} \Bigl( 1 + \frac{c_2}{m}
+ c_3 \Bigl( \frac{d}{m} \Bigr)^{\gamma-2} \Bigr) \\
&=& O(nd^{\gamma-1} m^{4-\gamma}).
\end{eqnarray*}
Hence,
\begin{eqnarray*}
\frac{\Var(D_2)}{\E[D_2]^2} &=& O \left( \frac{nd^{\gamma-1} m^{4-\gamma}}{n^2d^{2\gamma-2}m^{6-2\gamma}}\right)\\
&=&O \left( \frac{1}{nm}\left(\frac{m}{d}\right)^{\gamma-1}\right)
\end{eqnarray*}
Now, by (\ref{assumptions}), and the fact that $\gamma>2$ we have $(m/d)^{\gamma-1} = o(n)$. Since the maximum degree $m$ is assumed to  grow as a power of $n$, we have $\frac{\Var(D_2)}{\E[D_2]^2}=o(1)$. Note that $\E[D_2]\sim c nd^{\gamma-1}m^{3-\gamma}$, for $2<\gamma<3$.

Using Chebyshev's inequality, this establishes the claim of the lemma.
\hfill $\Box$

To complete the proof of Theorem \ref{thm:graphs}, it suffices to combine the results of lemma \ref{lem:total_deg_stats} and lemma \ref{lem:twostar_plg} with the fact that $\sum_v \pi_v^2=\frac{D_2+D}{D^2}$.

\section{Conclusion and future work}\label{se-conclusion}

In this work we have presented a simple model for the spread of epidemics where individuals are mobile. In this framework we were interested in the setting where there are two individuals one infected and one healthy both performing random walks on the network. Our preliminary investigation highlights the effect of  the topology on the spread of an epidemic, motivated by networking phenomena such as worms and viruses, failures, and dissemination of information. Under this natural
model, we provided an explicit relationship between the structure over which the walks are performed and the coincidence time of the two walkers. To this end we analysed both homogeneous (regular, complete and Erd\"{o}s-R\`enyi graphs) and heterogeneous (power-law graphs) networks. We pinpointed the existence of a phase transition for the coincidence time in the case of power-law networks depending on the parameter of the power-law degree distribution. We also derived bounds on the probability of infection.

As a final remark, we propose some several interesting directions to pursue the work presented
here. In our present model individuals are supposed to start their walks in stationary regime. This can be relaxed since the networks we study are expanders and thus random walks on such networks have nice mixing properties as illustrated in \cite{GaMaTow} through the computation of the isoperimetric constant of the underlying graphs. We also anticipate that similar results can be derived when considering $k$ walkers as long as $k$ is small with respect to $n$ the number of sites in the network.

\end{document}